\documentclass[a4paper,11pt]{article}
\usepackage[utf8]{inputenc}
\usepackage[dvipsnames]{xcolor}
\usepackage{amsmath,amsfonts,amssymb,amsthm} 
\usepackage{graphicx}
\usepackage{mathtools}
\usepackage[left=3cm,right=3cm,top=3.5cm,bottom=3.5cm]{geometry}
\usepackage[colorlinks=true]{hyperref}
\usepackage{mathptmx}
\usepackage{caption}
\usepackage[cal=cm,scr=euler,frak=lucida,frakscaled=1]{mathalfa}
\definecolor{Violet}{RGB}{238, 130, 238}
\definecolor{darkorange}{RGB}{255, 140, 0}
\newtheorem{theorem}{Theorem}[]
\newtheorem{proposition}{Proposition}[]

\newtheorem{corollary}{Corollary}

\newcommand{\IP}{\mathbb{P}}
\newcommand{\E}{\mathbb{E}}

\newcommand{\cT}{\mathcal{T}}

\renewcommand{\leq}{\leqslant}
\renewcommand{\geq}{\geqslant}
\newcommand{\1}{\text{\usefont{U}{bbold}{m}{n}1}}

\title{ \textsc{\bfseries  Parking on the Random Recursive Tree}}

\author{Alice Contat\thanks{LAGA - UMR CNRS 7539, Université Sorbonne Paris Nord, France and FSMP \hfill \href{mailto:alice.contat@math.cnrs.fr}{alice.contat@math.cnrs.fr}}\ , \; 
Lucile Laulin\thanks{Modal'X - UMR CNRS 9023, Université Paris Nanterre, 92000 Nanterre, France \hfill
\href{mailto:lucile.laulin@math.cnrs.fr}{lucile.laulin@math.cnrs.fr}}}
\begin{document}
\date{}
\maketitle

\abstract{We study the parking process on the random recursive tree. We first prove that although the random recursive tree has a non-degenerate Benjamini--Schramm limit, the phase transition for the parking process appears at density $0$. We then identify the critical window for appearance of a positive flux of cars with high probability. In the case of binary car arrivals, this happens at density $ \log (n)^{-2+o(1)}$ where $n$ is the size of the tree. This is the first work that studies the parking process on trees with possibly large degree vertices.}

\section{Introduction}
The parking process is now a well-studied algorithm in probability and combinatorics. It was introduced in the case of the line by Konheim and Weiss \cite{konheim1966occupancy} in the 1960's, motived by applications to hashing tables. More recently, it has been  generalized to the case of trees by Lackner and Panholzer \cite{LaP16}. Since then, the study of the parking process has focused  mainly on critical Bienaym\'e--Galton--Watson trees (e.g.\ \cite{chen2021parking,contat2020sharpness,contat2023these,contat2023parking,contat2021parking,curien2022phase,GP19}), on the complete binary tree \cite{aldous2022parking} or on possibly supercritical geometric trees \cite{chen2021enumeration,chen2024parking}. In all these cases, an interesting phase transition has been highlighted using a large range of techniques. 
However, this model has not been considered so far on trees which may have \textit{large degrees}. In this work, we start this new line of investigation by  focusing on the case of the (uniform) \emph{Random Recursive Tree} model. 

A recursive tree is a (non-plane) tree with vertices labeled from $1$ to $n$ and where the labels are increasing along the branches starting from the vertex labeled $1$ called the \emph{root}, see Figure \ref{fig:yuleRRT} (right). There are $(n-1)!$ recursive trees with $n$ vertices and we write $ \mathbb{T}_n$ for a uniform recursive tree of size $n \geq 1$. This model is of utmost importance in combinatorics and probability, and often serves as a benchmark for many algorithms, see e.g.~\cite{drmota2009random} and the references therein. It is easy to see that the degree of the root in $ \mathbb{T}_n$ is of order $\log n$, while the maximal degree of the vertices of $ \mathbb{T}_n$ is of order $ \log_2 (n)$, see \cite{devroye1995strong}. 

\paragraph{Parking on trees.} Recall the model of parking on a (fixed) rooted tree $ \mathbf{t}$, which represents the parking lot. The edges are thought of as one-way streets oriented towards the root. Random numbers of cars  arrive at the vertices of the tree and attempt to park. Each car prefers to park at its arrival vertex, but if that vertex is already occupied, it drives towards the root, following the oriented edges, until it encounters an empty vertex, where it parks.  If it reaches the root without being able to park, it exits the tree without parking. 

In the rest of the paper, we will always consider that the car arrivals, conditionally on the underlying tree, are i.i.d.\ on the vertices.  When their common distribution is $\mu$, we denote by $ \varphi (  \mathbf{t}, \mu)$ for the \emph{flux} of cars exiting $ \mathbf{t}$. In some cases, we also consider $ \psi ( \mathbf{t}, \mu) $ the number of cars visiting the root of $ \mathbf{t}$, so that $ \varphi ( \mathbf{t}, \mu) = ( \psi ( \mathbf{t}, \mu)-1 )_+ := \max (  \psi ( \mathbf{t}, \mu)-1,0)$. 
Our goal is to study the parking process  when the underlying tree $ \mathbf{t }$ has the law of $ \mathbb{T}_n$.

\paragraph{Phase transition.} For a broad class of tree and of car arrival distributions, the parking process undergoes a very rich phase transition. To fix ideas, we consider $ (\mu_ \alpha : \alpha \in [0,1])$ a stochastically increasing family of probability distributions of car arrivals, where the law $ \mu_ \alpha$ has mean $ \alpha  $  for all $ \alpha \in [0,1]$, and we assume that there exists a constant $K$ such that for all $ \alpha \in [0,1]$, the variance of $ \mu_ \alpha$ is smaller than $K$. To avoid trivialities, we also assume that as soon as $\alpha \in (0,1]$, there is a positive probability that two cars arrive on the same vertex (i.e.\ $ \mu_ \alpha \left( \{ 0,1\} \right) <1$), otherwise each car  parks on its arriving spot.  

We also consider a  family of (laws of) trees $ (T_n : n \geq 1)$ such that for all $n$, the tree $ T_n$ has size $n$. We assume that the sequence $(T_n : {n \geq 1})$ converges in the Benjamini--Schramm quenched sense (see Section \ref{sec:limloc} for a precise definition) towards a rooted infinite tree which has a single infinite spine almost surely. In this context, it has been shown \cite[Theorem 4.1]{contat2023these} that there exists a critical parameter $ \alpha_c \in [0,1]$ such that 

\begin{equation*} \frac{ \varphi (T_n, \mu_ \alpha)}{ n } \xrightarrow[n\to\infty]{( \mathbb{P})} \begin{cases}
	0 &\text{ if $\alpha < \alpha_c$}\\
	C_\alpha >0 &\text{ if $\alpha> \alpha_c$.}
\end{cases}
\end{equation*}
We say that the parking process is \emph{subcritical} (resp.\ \emph{supercritial}) when $ \alpha < \alpha_c$ (resp.\ \linebreak[4] {$\alpha > \alpha_c$}).
There are families of trees for which this critical parameter $ \alpha_c$ is explicit.  The first example that has been considered are Cayley trees \cite{contat2021parking,GP19,LaP16} with Poisson car arrivals, where the transition happens when the parameter of the Poisson distribution is $1/2$. It has been then extended to critical Bienaym\'e--Galton--Watson trees conditioned to have size $n$ \cite{contat2020sharpness,contat2023parking,curien2022phase} where the critical parameter $ \alpha_c$ only depends on the first two moments of the offspring and car arrival distributions. The characterization of the critical point is also known for supercritical Bienaym\'e--Galton--Watson tree with geometric distribution \cite{chen2024parking}, and for the infinite binary tree \cite{aldous2022parking} (even if this very last example does not fit  in this setting of converging sequence of finite trees).

\paragraph{``Trivial" phase transition.} Coming back to the context of random recursive trees, our first result is that the phase transition for the parking process on $ \mathbb{T}_n$ is ``trivial" in the sense that the critical parameter $ \alpha_c$ equals $0$. At first glance, one may think that this is an easy consequence of the fact the root vertex of $ \mathbb{T}_n$ has large degree (of order $ \log (n)$) as $ n$ goes to infinity. However, the phase transition involves the Benjamini--Schramm limit of the tree sequence, that is the local limit of the trees seen in the neighbourhood of a ``typical" uniform point, which is non-degenerate. For example, the degree of a vertex chosen uniformly at random in $ \mathbb{T}_n$ converges in law towards a random variable with a  Geometric distribution with parameter $1/2$.

\begin{theorem}\label{thm:macroflux}  
For any choice of family $(\mu_ \alpha : \alpha \in [0,1])$, the critical parameter for parking on the random recursive trees is zero. 
In other words, for any car arrival  law $ \mu$ such that $ \mu \left( \{ 0,1\} \right) <1$, there exists a positive constant $ C> 0$ such that 
$$ \frac{ \varphi( \mathbb{T}_n, \mu)}{n} \xrightarrow[n\to\infty]{ (\mathbb{P})} C.$$
\end{theorem}

\paragraph{Critical window for binary car arrivals.} 
The previous theorem says that the parking process is always supercritical. It is then natural to try to sophisticate the question and determine the order of magnitude of cars needed per vertex to see an outgoing flux of cars at the root with high probability. A way to do this is to allow the law of car arrivals to depend on~$n$: for each $n$, choose a parameter $ \alpha_n$ and study the sequence  $ (\varphi(  \mathbb{T}_n, \mu_{\alpha_n}) : n \geq 1)$. For the reader's convenience, let us first present our result in the case of binary car arrivals. That is we take for $ \alpha \in [0,1]$
\begin{equation*} \mu_{ \alpha}^{ \mathrm{Binary}} := \Big(1- \frac{\alpha}{2}\Big) \delta_0 + \frac{\alpha}{2} \delta_2.
\end{equation*}
In this case, the critical window is when $ \alpha_n$ is of order $ (\log(n))^{-2 + o(1)}$. 

\begin{theorem}\label{thm:binary}For all $ \delta >0$, we have 
\begin{eqnarray*} \varphi( \mathbb{T}_n, \mu_{ \alpha_n}^{ \mathrm{Binary}}) \xrightarrow[n\to\infty]{ ( \mathbb{P}) } \begin{cases}
	0 &\text{ if $\alpha_n \ll \log(n)^{-2}$}\\
	+ \infty  &\text{ if $\alpha_n \gg  \log(n)^{-2+ \delta}$}\end{cases}
\end{eqnarray*}
\end{theorem}
We used the notation $u_n \ll v_n$ if $u_n/v_n \to 0$ as $n \to \infty$. The random recursive trees $ \mathbb{T}_n$ are naturally coupled in an increasing fashion by letting the vertices arrived one by one, and connect the vertex labeled $n$ to a uniform vertex of $ \mathbb{T}_{n-1}$. This yields an increasing sequence $( \mathbb{T}_n : {  n \geq 1})$ of recursive trees, and we shall always work with them in what follows (and we keep the same notation).  If we assign the car arrivals to each vertex labeled $\{1,2, ...\}$ i.i.d.~of law $\mu$, this enables us to study the coupled version of the flux $\big(\varphi(\mathbb{T}_n,\mu) : n \geq 1\big)$ which are then increasing in $n$ almost surely. We can then define the time of appearance of the flux
 $ \theta (\mu, C )$ 
\begin{equation}\label{eq:deftime} \mbox{ for }C>0, \qquad \theta({\mu}, C) := \inf \Big\{ n \geq 0,  \varphi ( \mathbb{T}_n, \mu)  \geq C  \Big \}. \end{equation}
The former result can then be equivalently stated as:
\begin{corollary}\label{cor:binary} For all $C > 0$,  $$ \frac{\log \log \theta ( \mu_ \alpha^{ \mathrm{Binary}}, C )}{|\log \alpha|} \xrightarrow[ \alpha \to 0]{ ( \mathbb{P}) } \frac{1}{2}.$$
\end{corollary}

\paragraph{General critical window.} This phase  transition can be extend to more general context of stochastically increasing families of car arrival distributions. Again, we  consider a stochastically increasing family of probability measures $ (\mu_{ \alpha} : \alpha \in [0,1])$  representing the car arrivals with the same assumptions as above and two additional ones:
\begin{itemize}
\item the expectation of $ \mu_ \alpha$ is $ \alpha$; 
\item To avoid trivial cases, we assume that $ \mu_ \alpha ( \{ 0,1\}) <1$ for all $ \alpha \in ( 0,1]$.
\item the car arrivals are bounded, i.e.\ there exists $K > 0$ such that $\mu_ \alpha \left( \{ 0, \ldots, K \} \right) = 1$ for all $ \alpha \in [0,1]$;
\item For all $ 1 \leq k \leq K$, either $ \mu_ \alpha ( \{ k \}) = 0$ for all $ \alpha$ small enough (in which case we simply let $C_k :=0$), or there exists $C_k > 0 $ and $\beta_k \geq 1/k $ such that,  
\begin{equation}\label{eq:assumption_asymptotic} \mu_ \alpha ( \{ k \}) \underset{\alpha \to 0}{\sim} C_k  \alpha^{ \beta_k \cdot k}.\end{equation}
\end{itemize}
An important quantity, that we define now is 
\begin{equation*} \beta^* = \inf \{ \beta_k : 1 \leq k \leq K \mbox{ such that } C_k > 0\} .
\end{equation*}
Note that for the family $ \mu_{ \alpha}^{ \mathrm{Binary}}$, we have $ \beta^* = 1/2$. Note that $ \beta^* \leq 1$ since we must have  
\begin{equation*} \sum_{k = 1}^{K} k  C_k \alpha^{ \beta_k \cdot k} = \alpha (1 + o(1)), 
\end{equation*}
so that there exists at least one $k$ for which $\beta_k = 1/k$.
Our result is the following 
\begin{theorem}\label{thm:transition} Under the above assumptions, we have for all $ \delta > 0$, 
\begin{eqnarray*}  \varphi( \mathbb{T}_n, \mu_{ \alpha_n}^{})  \xrightarrow[n\to\infty]{ ( \mathbb{P})} \begin{cases}
	0 &\text{ if $\alpha_n \ll \log(n)^{-1/ \beta^* }$}\\
	+ \infty  &\text{ if $\alpha_n \gg  \log(n)^{-1/ \beta^* + \delta}$}\end{cases}
\end{eqnarray*}
As a consequence, we have for all $ C >0$, 
$$ \frac{ \log \log \left( \theta ( \mu_{ \alpha}, C) \right)}{ |\log ( \alpha) |} \xrightarrow[ \alpha \to 0]{ ( \mathbb{P})} \beta^*.$$
\end{theorem}
Note that Theorem \ref{thm:binary} and Corollary \ref{cor:binary} are a special case of Theorem \ref{thm:transition}, which means that it suffices to prove Theorem \ref{thm:transition}. The global strategy is to focus on the \emph{parked component} of the root i.e.\ the connected component of the root when we only keep the edges of the tree between pairs of vertices that both contain a parked car in the final configuration. On the one hand, when $ \alpha_n$ is too small, we use a first moment (inspired from the proof of  \cite[Theorem 3.4]{GP19}) to show that the expectation of the sum of the contributions of all possible configurations for the component of the root tends to $0$. On the other hand, we show that when $ \alpha_n$ converges to $0$ too slowly compared to $n$, then we can find inside the connected component of the root a multi-type Bienaym\'e--Galton--Watson tree. This tree is constructed by restricting to generations that are multiples of $k^* \in \mathrm{argmin} \{  \beta_k \}$ receiving $k^*$ cars on them. We then show that this structure, when reaching a certain height, induces a large flux of cars at the root.

\paragraph{Extensions.} Our results naturally  raise other questions. For example, one may wonder about car arrival distributions that are not bounded, such as Poisson or geometric distributions. We believe that Theorem \ref{thm:transition} remains valid in this case.  Our proof relies on the boundedness hypothesis only to establish the upper bound for the flux. 
 A natural next step would be to refine our results by determining the regime in which we obtain a finite but non-zero flux of outgoing cars. In fact, it would be interesting to determine the order of magnitude of $ \log \left( \theta( \mu_{ \alpha}, C)\right)$ or even that of $ \theta( \mu_{ \alpha}, C)$.

\paragraph{Acknowledgments.}We thank Christina Goldschmidt for an interesting discussion at an early stage of this project. This project benefits from the support of SMAI and CIRM as a selected ``BOUM Project 2023". We are grateful to Nicolas Curien for sharing us some figures of \cite{curien2024RWRG}.

\section{Trivial phase transition}

The goal of this section is to prove Theorem \ref{thm:macroflux}. The main ingredient is to use the characterization of the critical point given in \cite[Theorem 4.1]{contat2023these}, which only uses the Benjamini--Schramm limit of the tree sequence. To describe this limit, it will be more convenient to use the coupled random recursive trees $( \mathbb{T}_n : {n \geq 1})$  that is constructed from a (continuous) Yule tree process. Let us start by presenting this construction.

\subsection{Random recursive tree constructed from a Yule process} 

\begin{figure}[!h]
 \begin{center}
 \includegraphics[width=15cm]{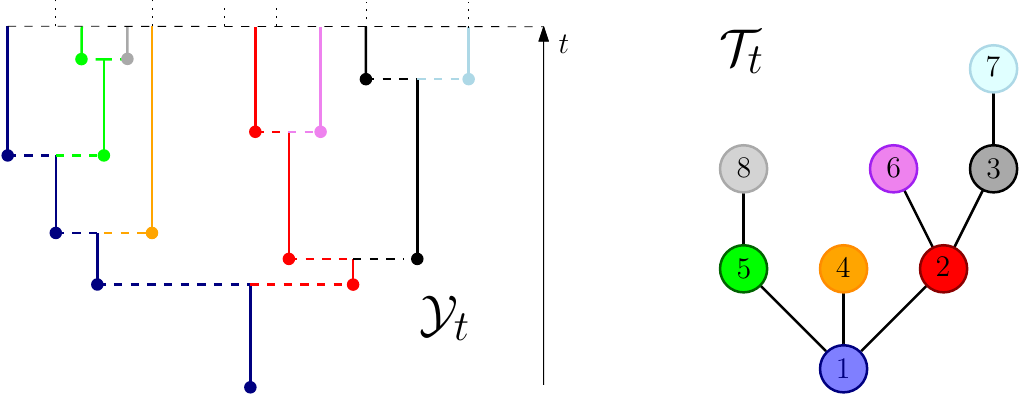}
 \caption{\label{fig:yuleRRT} On the left, an example of a Yule tree $ \mathcal{Y}_t$ cut at time $t$. On the right, the recursive $ \mathcal{T}_t$ tree constructed from this Yule tree. Each vertex is drawn with the same color as its corresponding branch in the Yule tree. }
 \end{center}
 \end{figure}

The starting point is a standard Yule tree.  The (infinite)  \emph{Yule tree}  is  a continuous-time branching process starting from one particle at time $0$ and each particle behaves independently of the other and lives an exponential time $ \mathcal{E} (1)$ of rate $1$ before dying and giving birth to $2$ new (independent) particles, see Figure \ref{fig:yuleRRT}. Formally, it can be represented as a \emph{plane} infinite binary tree where the edges are labeled with i.i.d.\ exponential random variables representing their lengths.

For each $ t \geq 0$, we denote by $ \mathcal{Y}_t$ the finite plane tree obtained by cutting the infinite Yule tree at time $ t$. 
From this Yule tree $\mathcal{Y}_t$, we can construct a recursive tree $ \mathcal{T}_t$ as follows:  we interpret all edges going to the left (recall that we have a plane orientation that enable us to distinguish left from right) as one vertex and at each branch point, the edge going to the right is the starting point of a new vertex, which we label with the smallest possible label and which is the children of the vertex corresponding to the edge left and down, see Figure~\ref{fig:yuleRRT}. Note in particular that the number of vertices of $ \mathcal{T}_t$ is the number of particles alive at time~$t$ in $ \mathcal{Y}_t$ and follows a geometric distribution with parameter $ \mathrm{e}^{-t}$ (see e.g.\ \cite{AN72}). Moreover, if $ 0 := \vartheta_1 < \vartheta_2 < \ldots < \vartheta_n < \ldots$ are the times of creation of particles in the Yule tree, then conditionally on $ ( \vartheta_i : {i \geq 0})$, we have 
\begin{equation}\label{eq:incr_coupling}\left( \mathcal{T}_{ \vartheta_n} : n \geq 1\right) \overset{(d)}{=} \left(\mathbb{T}_n: n \geq 1\right),\end{equation}
see e.g.\ \cite{holmgren2017fringe}. In particular, for all $ t \geq 0$, conditionally on its size $ | \mathcal{T}_t|$, the tree $ \mathcal{T}_t$ has the law of $ \mathbb{T}_{ | \mathcal{T}_t|}$.

\subsection{Description of the local limit of the RRT}\label{sec:limloc} 
We now describe the topology that we need in order to apply \cite[Theorem 4.1]{contat2023these}.
\paragraph{Benjamini--Schramm quenched topology.} Let $(T_n : n \geq 1)$ be a sequence of (possibly random) finite trees and $ (T_ \infty, \rho_ \infty)$ a (possibly infinite and random) rooted tree, that is a tree given with a distinguished vertex. For all $n \geq 1$, conditionally on $ T_n$, we let $ U_n$ and $V_n$ be two independent uniform vertices of $ T_n$. We say that the sequence $(T_n : n \geq 1)$ converges in the \emph{Benjamini--Schramm quenched} sense towards $( T_ \infty, \rho_ \infty)$ if for all $r \geq 0$ for all bounded functions $f$ and $g$ and such that the values of $ f ( \mathbf{g}, \rho)$ and $ g ( \mathbf{g}, \rho)$ only depend on the ball $ B_r ( \mathbf{g}, \rho)$ centered at $ \rho$ of radius $r$ in the graph~$ \mathbf{g}$,
\begin{equation}\label{eq:BSquenched} \mathbb{E} \left[ f( T_n, U_n) g(T_n, V_n) \right ] \xrightarrow[n\to\infty]{} \mathbb{E} \left[ f( T_ \infty, \rho_ \infty) \right ] \mathbb{E} \left[ g( T_ \infty, \rho_ \infty) \right ]. 
\end{equation}

\paragraph{Benjamini--Schramm quenched limit of the random recursive tree.} 
To apply \cite[Theorem 4.1]{contat2023these}, we have to know the Benjamini-Schramm quenched limit (if it exists) of the sequence of random recursive trees $( \mathbb{T}_n : n \geq 1)$. The existence of a Benjamini--Schramm limit (when we need only to consider one uniform point $ U_n$ and one function $ f$ instead of two) was already shown in \cite[Section 4]{aldous1991asymptotic}, which give us a natural candidate for the limit.  

\begin{figure}[!h]
 \begin{center}
 \includegraphics[width=\textwidth]{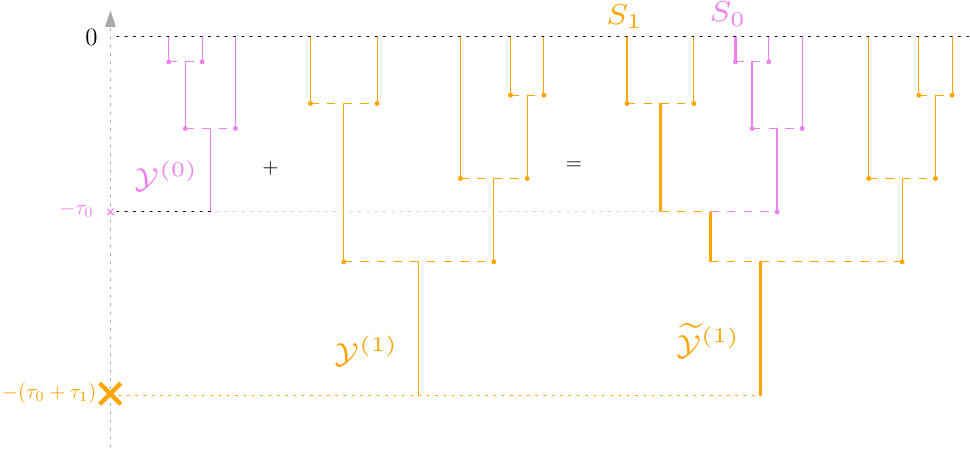}
 \caption{\label{fig:limloczoom} On the left side of the figure, we see the two Yule processes $\mathcal{Y}^{(0)}$ (in pink) and $\mathcal{Y}^{(1)}$ (in orange) that grow for respective times $\tau_0$ and $\tau_0+\tau_1$ in total. On the right side, the tree $ \widetilde{ \mathcal{Y}}^{(1)}$ obtain by gluing $\mathcal{Y}^{(1)}$ on the right-hand side of the left most branch of $ \widetilde{ \mathcal{Y}}^{(0)}= \mathcal{Y}^{(0)}$. The two thicker branches are the one that correspond to the vertex $ S_0$ (in pink) and $S_1$ (in orange).}
 \end{center}
 \end{figure}

\begin{figure}[!h]
 \begin{center}
 \includegraphics[width=7.5cm]{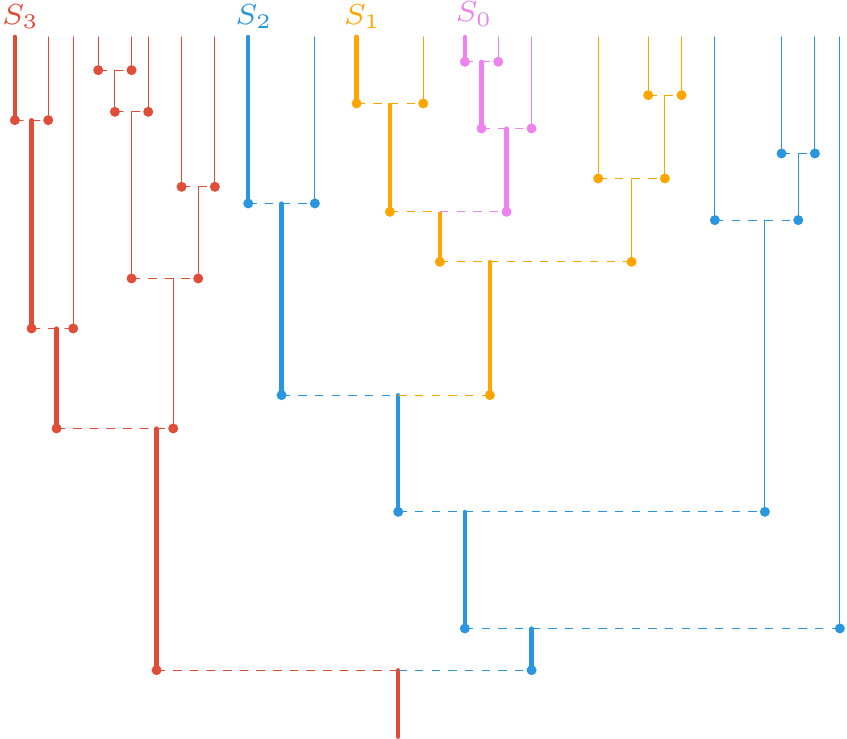} \quad  \includegraphics[width=7cm]{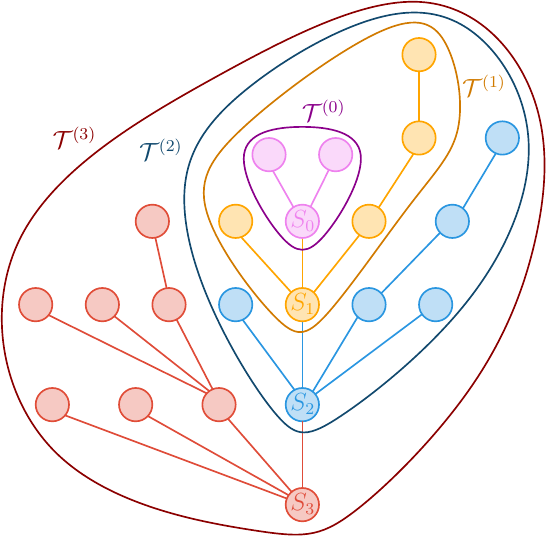}
 \caption{\label{fig:limloc} On the left side, the tree $ \widetilde{ \mathcal{Y}}^{(3)}$ obtained by gluing Yule tree as described above. The branch corresponding to the ancestors $S_0, S_1, S_2$ and $ S_3$ are drawn thicker. On the right, the recursive tree $ \mathcal{T}^{(3)}$ constructed from $ \widetilde{ \mathcal{Y}}^{(3)}$. We also highlight with colors the increasing (for the inclusion) sequence of recursive trees $ (\mathcal{T}^{(k)}: 0 \leq k \leq 3)$.
  }
 \end{center}
 \end{figure}

We thus denote by $ (\mathcal{T}_ \infty, S_0)$ the random rooted tree constructed as follows:
\begin{itemize}
\item Let $ (\tau_k)_{k \geq 0 }$ be a sequence of i.i.d.\ random variables which follow an exponential distribution $ \mathrm{Exp}(1)$.
\item Conditionally on $ (\tau_k)_{k \geq 0 }$, let $( \mathcal{Y}^{(k)} :{k \geq 0}) $ be a sequence of independent Yule tree processes such that the tree $\mathcal{Y}^{(k)}$ starts at time $ - \sum_{j = 0}^{k} \tau_k$ and stops its growth at time $0$.
\item Then we build recursively a sequence of trees $( \mathcal{T}^{(k)} : k \geq 0)$ from the sequence of Yule trees  $( \mathcal{Y}^{(k)} : {k \geq 0}) $ such that: the tree $ \mathcal{T}^{(0)}$ is the recursive tree constructed from the Yule tree $ \widetilde{ \mathcal{Y}}^{(0)}:= \mathcal{Y}^{(0)}$ (recall that it starts at time $ - \tau_0$ and  is stopped at time $ 0$). We call $S_0$ its root. For $ k \geq 1$, we see the vertex corresponding to the left-most branch of $\mathcal{Y}^{(k)}$ as the parent of $S_{k-1}$, the vertex corresponding to the left-most branch of  $ \widetilde{\mathcal{Y}}^{(k-1)}$. In other words, we graft the tree $\widetilde{\mathcal{Y}}^{(k-1)}$ on the right hand side of the left-most branch of $\mathcal{Y}^{(k)}$ to obtain $ \widetilde{\mathcal{Y}}^{(k)}$, see Figure \ref{fig:limloczoom}.  We then consider the recursive tree $\mathcal{T}^{(k)}$ built from $ \widetilde{\mathcal{Y}}^{(k)}$, rooted at a vertex which we denote by $ S_{k}$, see Figure \ref{fig:limloc}. 
\item The sequence of trees $( \mathcal{T}^{(k)} : k \geq 0)$ is increasing for the inclusion process, so that we can write $\mathcal{T}_ \infty$ for its limit. We see this tree as rooted at the vertex $ S_0$, which has an infinite line of ancestors $(S_k : {k \geq 0})$.
\end{itemize}
This way of describing $ \mathcal{T}_ \infty$ has the advantage of highlighting its genealogy, but there is another equivalent way to understand the construction of $ \mathcal{T}_ \infty$: Consider the semi-infinite branch $(- \infty, 0]$. Each side of this branch is decorated with a Poisson point process with intensity $1$, and at each atom appearing at time $- s$ with $s>0$, we start an independent Yule tree which stops it growth at time $0$. We can then contract all branches going to the left and get an infinite one-ended tree, whose distinguished vertex corresponds to the part of the infinite branch of the branch $(- \infty, 0]$ before the first atom on the left.

\begin{proposition}[] The tree $ \mathcal{T}_{ \infty}$ is the limit in Benjamini--Schramm quenched sense of $ \mathbb{T}_n$ as $n$ goes to infinity.
\end{proposition}

This has been shown in \cite[Example 6.1]{holmgren2017fringe}. Let us explain how to deduce \eqref{eq:BSquenched} from their results.
From \cite[Theorem 5.25]{holmgren2017fringe}, we know that for all  finite graph $ \mathbf{g}$ and for all $ r \geq 0$ we have 
\begin{eqnarray*}  \frac{1}{n}\sum_{v \in \mathbb{T}_n} \1_ { B_r( \mathbb{T}_n, v) = \mathbf{g}} \xrightarrow[n\to\infty]{( \mathbb{P})} \mathbb{P} \left( B_r( \mathcal{T}_{ \infty}, S_0) = \mathbf{g}\right).  \end{eqnarray*}
This implies, by Slutsky's theorem, that for all $r_1, r_2 \geq 0$ and graphs $ \mathbf{g}_1, \mathbf{g}_2$, we have 
$$ \frac{1}{n^2} \sum_{v_1,v_2 \in \mathbb{T}_n} \1_ { B_{r_1}( \mathbb{T}_n, v_1) = \mathbf{g}_1} \1_ { B_{r_2}( \mathbb{T}_n, v_2) = \mathbf{g}_2}  \xrightarrow[n\to\infty]{( \mathbb{P})} \mathbb{P} \left( B_{r_1}( \mathcal{T}_{ \infty}, S_0) = \mathbf{g}_1\right) \mathbb{P} \left( B_{r_2}( \mathcal{T}_{ \infty}, S_0) = \mathbf{g}_2\right).  $$
This implies \eqref{eq:BSquenched} by decomposing the functions $f$ and $g$ for all possible values of the pair $ ( B_r ( \mathbb{T}_n, v_1), B_r ( \mathbb{T}_n, v_2)) $ for $ v_1, v_2 \in \mathbb{T}_n.$

\subsection{Trivial phase transition $\alpha_c=0$}
The goal of this section is to show our Theorem	 \ref{thm:macroflux} stating that the phase transition for parking on the Random Recursive Tree is trivial. The rough idea is that, although a typical vertex of $ \mathbb{T}_n$ has a ``bounded degree", its $k$-th ancestor has degree of order $k$, and so has a probability of accommodating a car that is increasing in $k$. 

We will need to investigate the flux or number of cars driving through each vertex. Thus, for every tree $ \mathbf{t}$ and every vertex $v$ of $ \mathbf{t}$ and for all probability distribution $ \mu$, we denote by $\varphi_v ( \mathbf{t}, \mu)$ the number of cars going out from the vertex $v$ and $\psi_v ( \mathbf{t}, \mu)$ the number of cars passing by the vertex $v$, so that  $ \varphi_v ( \mathbf{t}, \mu)= (\psi_v ( \mathbf{t}, \mu)- 1)_+$. In particular, if $ \rho$ is the root of $ \mathbf{t}$, then we can omit the $\rho$ as an index in the notation such that  $\varphi ( \mathbf{t}, \mu) = \varphi_ \rho ( \mathbf{t}, \mu)$ is the flux of cars going out the tree and $\psi  ( \mathbf{t}, \mu) = \psi_ \rho ( \mathbf{t}, \mu)$ is the number of cars passing by the root. We recall that the \emph{clusters of parked cars} are the connected component of the tree when keeping only the edges between pairs of vertices which both contains a car in the final configuration.

 Let $ (\mu_ \alpha : \alpha \in [0,1])$ be a family of  probability distributions such that  $\mu_ \alpha$ has mean $ \alpha$ and such that for all $ \alpha \in (0,1]$, we introduce $ \delta( \alpha) :=  1 - \mu_ \alpha ( \{ 0,1\}) > 0$ . We are looking for
\begin{equation}
\label{def:alpha-crit}
	\alpha_c := \inf \left\{\alpha>0,\; \IP(\text{there exists an infinite size cluster of parked cars in $\cT_\infty$})>0\right\}.
\end{equation}
We know by \cite[Theorem 4.1]{contat2023these} that,
\begin{equation}
\label{eq:flux}
\frac{\varphi( \mathbb{T}_n)}{n} \overset{( \mathbb{P})}{\underset{n\to\infty}\longrightarrow} 
\begin{cases}
	0 &\text{ if $\alpha < \alpha_c$}\\
	C_\alpha >0 &\text{ if $\alpha> \alpha_c$}
\end{cases}
\end{equation}
 We now have now all the ingredients to prove Theorem \ref{thm:macroflux}, i.e.\ that  $ \alpha_c = 0$. 
\begin{proof}  Let us fix $ \alpha \in (0,1]$. We want to show that the probability that there exists an infinite cluster of parked cars in $ \mathcal{T}_ \infty$ is positive. Let us recall from Section \ref{sec:limloc} the notation and construction of the infinite tree $ \mathcal{T}_{ \infty}$ as the limit for the inclusion process of the trees $ (\mathcal{T}^{(k)} : k \geq 0) $. Note in particular that conditionally on the $ (\tau_k)_{k \geq 0}$, the trees $ \mathcal{T}^{(k)} \setminus \mathcal{T}^{(k-1)}$ are independent and have law $ \mathcal{T}_{\sum_{i=0}^k \tau_i}$.  
We thus start by considering a fixed $t>0$ and give an upper bound for the probability that the root vertex $ \rho$ of a recursive tree $ (\mathcal{T}_t, \rho)$ constructed from a Yule tree $ \mathcal{Y}_t$ remains a free spot in the final configuration of cars. Actually, we will only use  the children of the root:  if the root remains a free spot, then each child of the root must have at most $1$ car arriving on it.
Recalling that the degree of the root, which we denote by $ \mathrm{deg} (\rho)$, follows a Poisson distribution with parameter $t$, we obtain 
\begin{align*}
	\mathbb{P}( \psi ( \mathcal{T}_t, \mu_ \alpha) = 0) &\leq \mathbb{E} \left[ \prod_{i= 1}^{ \mathrm{deg}(\rho)} \1_{\mbox{at most $1$ car is arriving on the $i$-th children of $ \rho$ in $ \mathcal{T}_t$}} \right ] \\
	&=  \E\left[(1- \delta(\alpha))^{\deg(\rho)}\right] \\
	&= \E\left[\sum_{k\geq 0}(1- \delta(\alpha))^{k} \frac{t^{k}}{k!}\mathrm{e}^{-t}\right] \\
	&= \mathrm{e}^{-t}\mathrm{e}^{(1- \delta(\alpha))t} = \mathrm{e}^{- \delta(\alpha) t}
\end{align*}
We now use this to  give an upper bound on the probability that $S_k$ is a free spot in $ \mathcal{T}_{ \infty}$. Note that it only depends on $ \mathcal{T}^{(k)}$ since $ \mathcal{T}^{(k)}$ contains all the vertices of $ \mathcal{T}_{ \infty}$ that are above $ S_k$. Recall that the $\tau_j$ have an exponential distribution with parameter $1$. So 
\begin{eqnarray*}
\mathbb{P} \left( \psi_{S_k} ( \mathcal{T}_{ \infty}, \mu_ \alpha) = 0\right) &\leq& \mathbb{E} \left[ \mathbb{E} \left[ \1_{ \psi_{S_k} ( \mathcal{T}^{ (k)}, \mu_ \alpha) = 0} \big| (\tau_j)_{0 \leq j \leq k} \right ] \right ] \\
&\leq& \mathbb{E} \left[ \mathbb{E} \left[  \mathrm{e}^{-  \delta (\alpha) \cdot \sum_{j = 0}^k \tau_j} \big| (\tau_j)_{0 \leq j \leq k} \right ] \right ] \\
&=& \int_{( \mathbb{R}_{ \geq 0})^{k+1}}\mathrm{e}^{- \delta(\alpha) \cdot \sum_{j = 0}^k t_j}  \prod_{j = 0}^{k } \mathrm{d}t_j \cdot \mathrm{e}^{- t_j}\\
&=& \int_{( \mathbb{R}_{ \geq 0})^{k+1}} \prod_{j = 0}^{k } \mathrm{d}t_j \cdot \mathrm{e}^{- (1+ \delta(\alpha)) t_j} \\
&=& \left( \frac{1}{1+ \delta(\alpha)}\right)^k
\end{eqnarray*}

In particular, we obtain from this computation that 
\begin{equation*}
	\sum_{k\geq1} \IP(\psi_{S_k} ( \mathcal{T}_{ \infty}, \mu) = 0)<\infty.
\end{equation*}
Consequently, we can conclude using Borel-Cantelli Lemma's that
\begin{equation*}
	\mathbb{P}\big(\nexists\text{ cluster of infinite size in } \cT_\infty \big) = \mathbb{P}(\forall n,\ \exists k\geq n,\ \psi_{S_k} ( \mathcal{T}_{ \infty}, \mu) = 0\big)= 0.
\end{equation*}
\end{proof}

\section{Upper bound for the flux}
We now want to refine the transition and study the behaviour of the flux of outgoing cars when  both $ n$ goes to $+ \infty$ and $ \alpha$ goes to $ 0$. Recall our assumptions and in particular the asymptotic~\eqref{eq:assumption_asymptotic} on the law $ \mu_ \alpha$ as $ \alpha$ goes to $0$.
The goal of this section is to provide an upper bound for the expectation of flux and for $ \theta({ \mu_ \alpha}, K)$ for $K > 0$. 
We start by looking at the continuous time version $ (\mathcal{T}_t : {t \geq 0})$  of the recursive tree and we will then come back to the discrete time version $ (\mathbb{T}_n : n \geq 1)$.
We prove the following proposition, from which we will deduce the first part of Theorem \ref{thm:transition}.
\begin{proposition}[subcritical regime]\label{prop:subcritic}There exists a constant $c>0$ such that if $(\alpha_t : t  \geq 0 )$ is such that for all $t$ large enough $\alpha_t^{\beta^*}  t \leq c$, then 
$$ \mathbb{E} [ \varphi ( \mathcal{T}_t, \mu_{\alpha_t})] \xrightarrow[t\to\infty]{} 0.$$
\end{proposition}

\subsection{Fully parked trees and decomposition}

\begin{figure}[!h]
 \begin{center}
 \includegraphics[width=11cm]{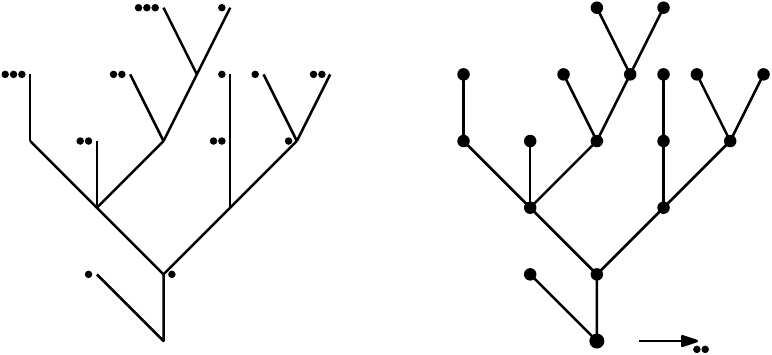}
 \caption{\label{fig:FPT}Example of a fully parked tree with $18$ vertices and $20$ cars arriving on it: on the left, the car arriving configuration and on the right, the final configuration where all spots are occupied and $2$ cars are going out from the tree and contributing to the flux.}
 \end{center}
 \end{figure}

To give an upper bound for the expectation flux of outgoing cars, we introduce the notion of \emph{fully parked trees}. A fully parked tree of size $n$ with $m$ cars is a plane tree with $n$ vertices carrying $m \geq n$ cars such that in the final configuration, all vertices contain a parked car, see Figure \ref{fig:FPT}. We denote by $ \mathrm{FPT}(n,m)$ the set of fully parked tree with $n$ vertices and $m$ cars. Note that, for a fully parked tree of size $n$ with $m$ cars, there is $m-n +1$ cars passing by the root.

An important remark is that all cars that pass by the root are arriving in the fully parked cluster of the root. Thus, we can decompose the flux according to the shape and car configuration of the cluster of the root.

\subsection{Proof of Proposition \ref{prop:subcritic}}
We first consider $ \alpha$ small enough so that for all $ 1 \leq k \leq K$, we have 
\begin{equation}\label{eq:alpha_small} \mu_{ \alpha} (k) \leq 2 C_k \alpha^{ \beta_k \cdot k} .
\end{equation}

Let us fix a fully parked tree $ \mathbf{t}$ in $ \mathrm{FPT}(n,m)$ for some $m \geq n \geq1$ and if $ v$ is a vertex of $ \mathbf{t}$, we write $ c_v$ for the number of cars arriving at $v$. Note that the sum of the $c_v$'s is $m$.  We want to give an upper bound for the probability that $ \mathbf{t}$ (as a plane tree with a car arrival decoration) is contained the parked connected component $ \mathcal{C}( \rho, \mathcal{T}_t )$ of the root in $ \mathcal{T}_t$. Note that the tree $ \mathcal{T}_t$ is a priori not planar, but it has a natural planar ordering coming from the labels of the vertices (we order the children from left to right by increasing label). Recall that $ \mathcal{T}_t$ is constructed from a Yule tree $ \mathcal{Y}_t$. For each $ v $ vertex of $ \mathbf{t}$ and $ 1 \leq i \leq c_v$, we integrate over all possible values of $ t_{v}^{i}$ the time corresponding to the creation of the $i$-th children of $v$ in $\mathbf{t}$ in $ \mathcal{T}_t$.

Given the family of $(t_v^i : v \in \mathbf{t}, 1 \leq i \leq c_v)$, then the probability that this embedded tree is contained in $ \mathcal{C} ( \rho, \mathcal{T})$, the parked component of the root of $ \mathcal{T}_t$, can be bounded above by the product of the probability of its car arrivals.  In particular, since $ \alpha$ is small enough so that Equation \eqref{eq:alpha_small} is satisfied, we get 
\begin{eqnarray*} \prod_{v \in \mathbf{t}} \mu (c_v) &\leq& \prod_{v \in \mathbf{t}} \left(2\max_{1 \leq j \leq K} C_j\right) \alpha^{\beta_{c_v} c_v} \\
&\leq& \left(2\max_{1 \leq j \leq K} C_j\right)^n \alpha^{ \sum_{v \in \mathbf{t}} \beta_{c_v} c_v} \\
&\leq& \left(2\max_{1 \leq j \leq K} C_j\right)^n \alpha^{\beta^* \sum_{v \in \mathbf{t}}  c_v} \\
&\leq& \left(2\max_{1 \leq j \leq K} C_j\right)^n \alpha^{\beta^* m}.
\end{eqnarray*}
Integrating over the $(t_v^i : v \in \mathbf{t}, 1 \leq i \leq c_v)$, we obtain
\begin{eqnarray*} \mathbb{P} \left( \mathcal{C} ( \rho, \mathcal{T}_t) \supseteq \mathbf{t}\right) &= &\mathbb{E} \left[ \1_ {\mathcal{C} ( \rho, \mathcal{T}_t)\supseteq \mathbf{t} }\right ] \\
&\leq& \prod_{v \in \mathbf{t}} \prod_{1 \leq i \leq c_v} \int_{0}^{t} \mathrm{d} t_{v}^i \left(2\max_{1 \leq j \leq K} C_j\right)^n \alpha^{\beta^* m} \\
&=& t^{n-1} \left(2\max_{1 \leq j \leq K} C_j\right)^n \alpha^{\beta^* m} 
\end{eqnarray*}

Now, we come back to  the expectation of $ \psi (\mathcal{T}_t, \mu_ \alpha)$. By summing over all fully parked trees that can be contained in the cluster of the root (which is itself a parked tree), we obtain
 \begin{eqnarray*} \mathbb{E} \left[ \psi (\mathcal{T}_t, \mu_ \alpha) \right ]  &\leq& \sum_{n \geq 0} \sum_{ m \geq n} (m-n+1) \sum_{ \mathbf{t} \in \mathrm{FPT} (n,m)} \mathbb{P} \left( \mathcal{C}( \rho, \mathcal{T}_t) \supseteq \mathbf{t}\right) \\
 & \leq&\sum_{n \geq 0} \sum_{ m \geq n} (m-n+1) \sum_{ \mathbf{t} \in \mathrm{FPT} (n,m)} t^{n-1} \left(2\max_{1 \leq j \leq K} C_j\right)^n \alpha^{\beta^* m}
 \end{eqnarray*}

Moreover, the number of possible shapes for the tree $ \mathbf{t}$ in $ \mathrm{FPT}(n,m)$ is the number of plane trees with $n$ vertices, which can be bounded by $ \mathrm{Cst} 4^n$ for some $ \mathrm{Cst}>0$, and there are at most $ (K+1)^n$ possible configurations of car arrivals since  the number of cars arriving at a specific vertex is bounded by $K$. Thus the cardinal of $ \mathrm{FPT}(n,m)$ is at most $ \mathrm{Cst} (4(K+1))^n$ and we write $C := 8 (K+1)\max_{1 \leq j \leq K} C_j$. 
Now we take a sequence $ ( \alpha_t : t \geq 0)$ such that  $   t \alpha_t ^{  \beta^*} \leq (1/ 2C)$ for all $t$ large enough. Then, for all $t$ large enough, we obtain
 \begin{eqnarray*} \mathbb{E} \left[ \psi (\mathcal{T}_t, \mu_ {\alpha_t}) \right ]
 &  \leq&\sum_{n \geq 0} \sum_{ m \geq n} \frac{\mathrm{Cst}}{t} (m-n+1) (C t)^n \alpha_t^{\beta^* (m)} \\
 &  \leq&\sum_{n \geq 0}  \frac{\mathrm{Cst}}{t}  \left(C t \right)^n  \sum_{ m \geq n}  (m+1) \alpha_t^{\beta^* m} \\
  &    \leq& \frac{\mathrm{Cst}}{t (1 - \alpha_t^{\beta^*})^2}\sum_{n \geq 0}   (n+1) \left( Ct \alpha_t^{ \beta^*}\right)^n \\
 &    \leq& \frac{\mathrm{Cst}}{t (1 - \alpha_t^{\beta^*})^2} \cdot \frac{1}{(1 - C  t \alpha_t^{  \beta^*})^2 } 
 \end{eqnarray*}
This proves the proposition with $c = 1/(2C)$.

\section{Lower bound for the flux}

We now want to show that if $ \alpha_n$ goes to $0$ too slowly, then the flux converges towards $ +\infty$ in probability. Again, we start by considering $ \mathcal{T}_t$ instead of $ \mathbb{T}_n$.

Recall the assumption on \eqref{eq:assumption_asymptotic} on the family $ (\mu_ \alpha: \alpha \in [0,1])$ as $ \alpha \to 0$.
Let $ \gamma> 0$ and $c> 0$ be such that for all $ 0<\alpha<1$, 
\begin{equation}\label{eq:gamma} 1 - \mu_ \alpha ( \{ 0,1\}) \geq c\alpha^{ \gamma}.\end{equation}
Moreover, we let $ k^* \in \mathrm{argmin} \{ \beta_k\} $, so that $ \beta^* = \beta_{k^*}$.  Let $ 0< \delta < 1$ and $ (\alpha_t:  t \geq 0)$ be a sequence such that for all $t$ large enough, we have \begin{equation}\label{eq:alpha_large} \alpha_t^{ \beta^*} t \geq t^{ \delta / k^*}.\end{equation} The global strategy is to construct a large tree made of vertices, which are at generations that are multiples of $k^*$ and that have $k^*$ cars arriving on it. Then when it reaches a certain height, we can find children of the vertices of this tree which will produce a large flux of cars.

\begin{figure}[!h]
 \begin{center}
 \includegraphics[height=7.2cm]{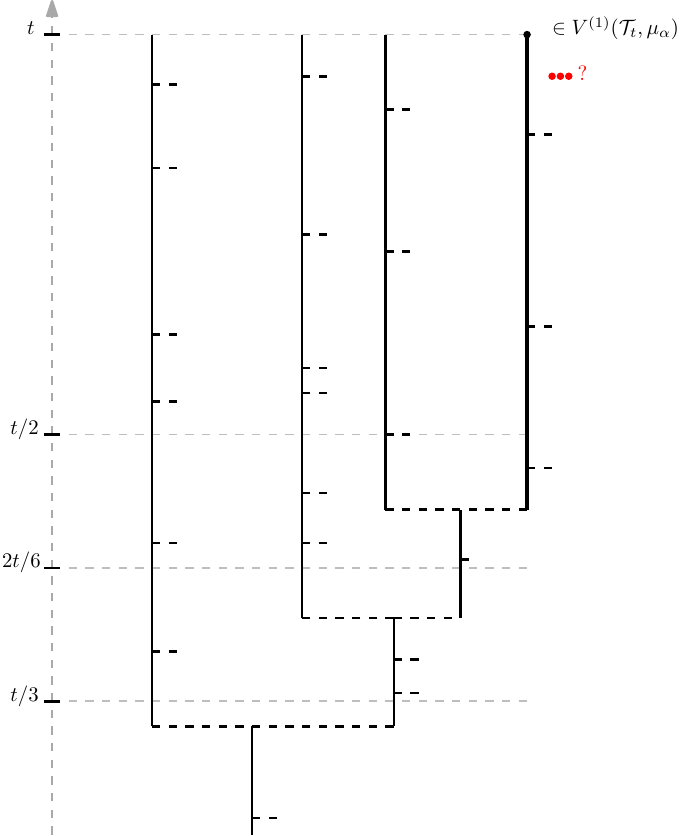} \hspace{0.2cm} \includegraphics[height=7.2cm]{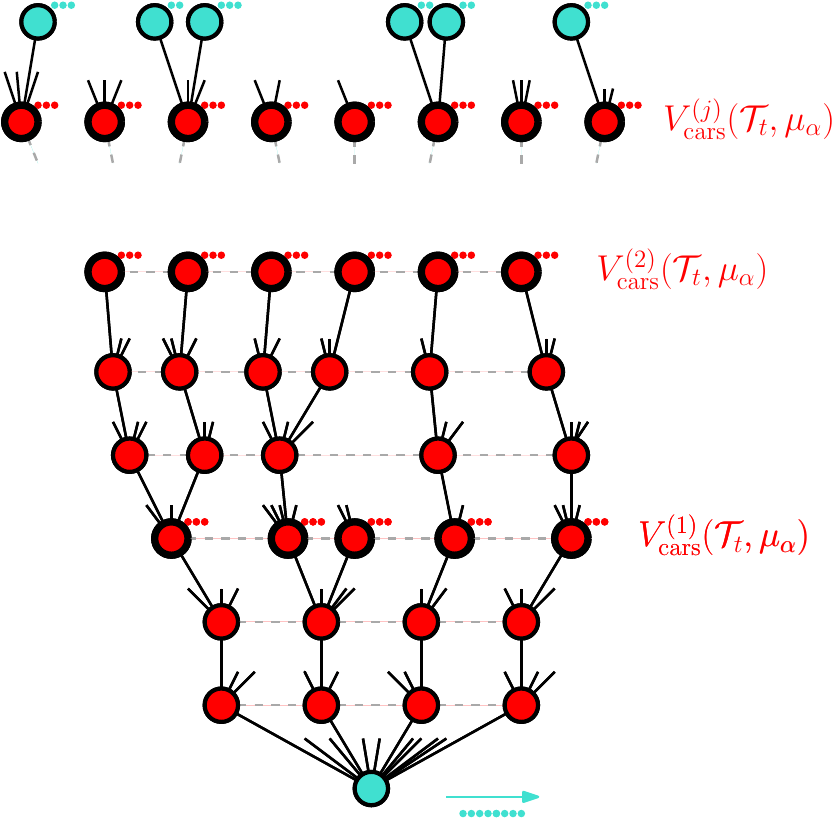}
 \caption{\label{fig:verticesV} On the left, example of a spliting in the Yule tree that creates a vertex in $ V^{(1)}( \mathcal{T}_t, \mu_ {\alpha_t})$. On the right, example of the construction of the sets $ ( V^{(\ell)}_{ \mathrm{cars}}( \mathcal{T}_t, \mu_ {\alpha_t}) : \ell \geq 1) $ when $ k^* = 3$. In our example, there are $6$ children of  $V^{(j)}_{ \mathrm{cars}}( \mathcal{T}_t, \mu_ {\alpha_t}$ on which at least $2$ cars arrives (with $15$ cars on them in total, so that there are at least $15-6=9$ cars visiting the root and $8$ outgoing cars.}
 \end{center}
 \end{figure}

For $t > 0$, let us denote by $ V^{(1)}( \mathcal{T}_t)$ the vertices of $ \mathcal{T}_t$ which are at height $k^*$ and which arrived before time $ t/2 $ in $ \mathcal{Y}_t$ and  by $V^{(1)}_{ \mathrm{cars}}( \mathcal{T}_t, \mu_ {\alpha_t})$ the vertices of $V^{(1)}( \mathcal{T}_t)$ that have $ k^*$ cars arriving on it, see Figure \ref{fig:verticesV}. Note that all vertices that are ancestors of a vertex of $V^{(1)}_{ \mathrm{cars}}( \mathcal{T}_t, \mu_ {\alpha_t})$ (except possibly the root) contain by construction a parked car in the final configuration, even if we use the cars arriving on vertices which are at height at most $ k^*$. 
We lower bound the size of $ V^{(1)}( \mathcal{T}_t)$ by the number of line of $k^*$ successive descendants of the root such that for $1 \leq i \leq k^*$, the $i$-th descendant of the root arrives at most at time $  t/(2k^*)$ after the $(i-1)$-th descendant. 
Hence, the subtree of $ \mathcal{T}_t$ spanned by $ V^{(1)}( \mathcal{T}_t)$ is stochastically included in a Bienaym\'e--Galton--Watson tree with offspring distribution $ \mathrm{Poisson} (t/(2k^*))$, so the expected number of descendants at height $k^*$ is $(t/(2k^*))^{k^*}$. 
By successive applications of the law of large numbers, we have 
$$ \mathbb{P} \left( |V^{(1)}( \mathcal{T}_t)| \geq \frac{1}{2}\left( \frac{t}{2k^*}\right)^{k^*}\right) \xrightarrow[t\to\infty]{} 1.$$
Among the vertices of  $V^{(1)}( \mathcal{T}_t)$, each of them contains $k^*$ cars with probability $ \mu_ {\alpha_t} ( k^*)$ independently.
Thus, 
$$\mathbb{P} \left( |V^{(1)}_{ \mathrm{cars}}( \mathcal{T}_t, \mu_ {\alpha_t})| \geq  \frac{1}{4}\left( \frac{t}{2 k^*}\right)^{k^*} \cdot \mu_{\alpha_t} (k^*)\right) \xrightarrow[t\to\infty]{} 1.$$
In particular, by our choice of $t$ and $ \alpha_t$, there exists $ M > 0$ such that for all $ \varepsilon >0$, for all $t$ large enough, we have
$$ \mathbb{P} \left(| V^{(1)}_{ \mathrm{cars}}( \mathcal{T}_t, \mu_ {\alpha_t}) |\geq M t^{ \delta  } \right)  \geq 1 -  \varepsilon.$$

Now, for $ \ell \geq 1$, we define recursively $  V^{(\ell+1)}( \mathcal{T}_t, \mu_ {\alpha_t})$ as the vertices that are at height $k^*$ above a vertex of $V^{(\ell)}_{ \mathrm{cars}}( \mathcal{T}_t, \mu_ {\alpha_t})$ and that arrive before time $ t - t/2^{\ell} $ in $ \mathcal{Y}_t$ and $V^{(\ell+1)}_{ \mathrm{cars}}( \mathcal{T}_t, \mu_ {\alpha_t})$ the vertices of $V^{(\ell+1)}( \mathcal{T}_t, \mu_ {\alpha_t})$ that have $ k^*$ cars arriving on it. 

Recall that $ \gamma$ is defined so that the family $( \mu_ \alpha : \alpha \in [0,1])$  satisfies \eqref{eq:gamma}. Let $j$ be an integer such that $ j >   \gamma ( k^* - \delta)/ ( k^* \beta^* \delta)$, so that, by \eqref{eq:alpha_large}, we have $t^{ \delta j } \alpha_t^{ \gamma} > 1$ for $t \geq 1$. Again by successive applications of the law of large number, we can find $ m \geq 0$ such that for all $ \varepsilon > 0$, for all $t$ large enough,  we have  
$$ \mathbb{P} \left(| V^{(j)}_{ \mathrm{cars}}( \mathcal{T}_t, \mu_ {\alpha_t}) |\geq m t^{ j \delta } \right)  \geq 1 - 2 \varepsilon.$$

Lastly, look at the \emph{children} of the vertices of $V^{(j)}_{ \mathrm{cars}}( \mathcal{T}_t, \mu_ {\alpha_t})$. The number of children of the vertices of $V^{(j)}_{ \mathrm{cars}}( \mathcal{T}_t, \mu_ {\alpha_t})$ is stochastically larger than the sum of $|V^{(j)}_{ \mathrm{cars}}( \mathcal{T}_t, \mu_ {\alpha_t})|$ i.i.d.\ Poisson random variables with parameter $ t/2^j$. Hence with probability at least $1- \varepsilon$, there are at least $ m c t^{ \delta j}{\alpha_t}^{ \gamma}$ children of the vertices of $V^{(j)}_{ \mathrm{cars}}( \mathcal{T}_t, \mu_{\alpha_t})$ that receive at least two cars. By construction, the ancestors of the vertices of $V^{(j)}_{ \mathrm{cars}}( \mathcal{T}_t, \mu_ {\alpha_t})$ (except possibly the root) are occupied by cars that arrives on vertices which are at most at height $ j k^*$. Thus, all children  of vertices of $V^{(j)}_{ \mathrm{cars}}( \mathcal{T}_t, \mu_{\alpha_t})$ that receive at least two cars will contribute to $ \psi ( \mathcal{T}_t, {\alpha_t})$. 
In particular, for all $t$ large enough,
$$\mathbb{P} \left( \psi ( \mathcal{T}_t, {\alpha_t}) \geq \frac{m c}{2^{j+1}}   t^{ \delta j  +1} \alpha_t^{ \gamma} \right) \geq 1 - 4 \varepsilon.$$ 
Now let $ C \geq 0$. Then, for all $t $ large enough, we have
\begin{eqnarray} \mathbb{P} \left(\psi ( \mathcal{T}_t, \alpha_t) \geq C \right) & \geq&  \mathbb{P} \left(\psi ( \mathcal{T}_t, \alpha_t) \geq \frac{ m c}{2^{j+1}} t  \right)\notag \\
&\geq&\mathbb{P} \left(\psi ( \mathcal{T}_t, \alpha_t) \geq \frac{m c}{2^{j+1}}   t^{ \delta j  +1} \alpha_t^{ \gamma} \right)\notag \\
&\geq& 1 - 4 \varepsilon. \label{eq:inf}
\end{eqnarray}
Equivalently, we have $$ \psi ( \mathcal{T}_t, \alpha_t) \xrightarrow[t\to\infty]{( \mathbb{P})} + \infty.$$

\section{Coming back to the discrete}

It only remains to come back to the discrete sequence of recursive trees $( \mathbb{T}_n : n \geq 1)$.
\proof[Proof of Theorem \ref{thm:transition}] 
Recall that the sequence $ ( \mathbb{T}_n : n \geq 1)$ is constructed so that it is increasing for the inclusion. Thus, for every fixed $ \alpha \in [0,1]$, the sequence $ (\psi( \mathbb{T}_n, \mu_ \alpha) : n \geq 1)$ is non-decreasing. 
In particular, recalling \eqref{eq:incr_coupling}, we have for all $ t \geq 0$ and $n \geq 0$, 
\begin{eqnarray*} \mathbb{E} \left[ \psi ( \mathcal{T}_t,  \mu_ \alpha) \right ] &\geq& \mathbb{E} \left[ \psi ( \mathcal{T}_t,  \mu_ \alpha) \1_{|\mathcal{T}_t| \geq n }\right ] \\
&\geq& \mathbb{E} \left[ \psi ( \mathbb{T}_n,  \mu_ \alpha) \1_{|\mathcal{T}_t| \geq n }\right ] \\
&\geq& \mathbb{E} \left[ \psi ( \mathbb{T}_n,  \mu_ \alpha) \right ]  \mathbb{P} \left( {|\mathcal{T}_t| \geq n } \right).
\end{eqnarray*}
Recall that the size $|\mathcal{T}_t| $  of $ \mathcal{T}_ t$ follows a geometric distribution with parameter $ \mathrm{e}^{-t}$. Thus, 
\begin{eqnarray*} \mathbb{P} \left( |\mathcal{T}_t| \geq n \right) = \left( 1 - \mathrm{e}^{- t}\right)^{n-1}
\end{eqnarray*}
Hence, taking $t = 2 \log n$, we obtain 
\begin{eqnarray*} \mathbb{P} \left(\psi ( \mathbb{T}_n,  \mu_ \alpha) > 0 \right) &=& \mathbb{P} \left(\psi ( \mathbb{T}_n,  \mu_ \alpha) \geq 1 \right)\\
& \leq&  \mathbb{E} \left[ \psi ( \mathbb{T}_n,  \mu_ \alpha) \right ] \\
&\leq& \frac{\mathbb{E} \left[ \psi ( \mathcal{T}_t,  \mu_ \alpha) \right ] }{ \mathbb{P} \left( {|\mathcal{T}_t| \geq n } \right)} \\
&\leq& \mathbb{E} \left[ \psi ( \mathcal{T}_t,  \mu_ \alpha) \right ]  \cdot   \left( 1 - \frac{1}{n^2}\right)^{-(n-1)}
\end{eqnarray*}
Thus, if $ \alpha_n \leq (2 \log (n)/ c)^{ - 1/ \beta^*}$, we have $ \alpha_n^{\beta^* } t  \leq c$, so that we can apply Proposition \ref{prop:subcritic} and the right hand-side goes to $0$ as $n$ goes to $ + \infty$. This concludes the proof of the first part of Theorem \ref{thm:transition}.
It also implies that for all $ C > 0$, 
$$ \liminf_{ \alpha \to 0} \frac{ \log \log \theta ( \mu_ \alpha, C) }{ | \log \alpha|} \geq \frac{1}{2}. $$

Now for the upper bound, we let  $ ( \alpha_n: n \geq 0)$ be such that  $ \alpha_n \geq \log (n)^{- \frac{1}{ \beta^*} (1 - \delta/ k^*)} $ for all $ n$ large enough. Then, using again the fact that for all $ \mu$, the flux $(\varphi( \mathbb{T}_n, \mu) : n \geq 1)$ is non-decreasing, we have for all $ n \geq 1$,

\begin{eqnarray*} \mathbb{P} \left( \varphi( \mathbb{T}_n, \mu_ {\alpha_n }) \geq K \right) &\geq& \frac{ \mathbb{P} \left( \varphi( \mathcal{T}_{\log(n)}, \mu_{\alpha_n} ) \geq K\ \&\ | \mathcal{T}_{\log(n)}| \leq n \right)}{ \mathbb{P} \left( | \mathcal{T}_{\log(n)} | \leq n \right)} \\
& \geq& 1 - \frac{1 -\mathbb{P} \left( \varphi( \mathcal{T}_{\log(n)}, \mu_ {\alpha_n}) \geq K \right) }{ \mathbb{P} \left( | \mathcal{T}_{\log(n)} | \leq n \right)}
\end{eqnarray*} 
Let $ \varepsilon> 0$. By Equation \eqref{eq:inf}, for all $n$ large enough, we have 
$$1 -\mathbb{P} \left( \varphi( \mathcal{T}_{\log(n)}, \mu_ {\alpha_n}) \geq K \right) \leq \varepsilon.$$
Moreover, if $n$ is large enough, we have  \begin{eqnarray*} \mathbb{P} \left( | \mathcal{T}_t | \leq n \right) = 1 - \left(1- \frac{1}{n }\right)^n \geq \frac{1}{2},
\end{eqnarray*}
which gives 
\begin{eqnarray*} \mathbb{P} \left( \varphi( \mathbb{T}_n, \mu_ {\alpha_n }) \geq K \right) & \geq& 1 - 2 \varepsilon,\qquad \mbox{and}\qquad  \varphi( \mathbb{T}_n, \mu_ {\alpha_n }) \xrightarrow[n\to\infty]{( \mathbb{P})} + \infty .\end{eqnarray*}
It also implies that 
$$ \limsup_{ \alpha \to 0} \frac{ \log \log  \theta ( \mu_ \alpha, C)}{ | \log \alpha|} \leq \frac{1}{2},$$
for all $ C > 0$, which concludes the proof of Theorem \ref{thm:transition}.

\endproof

\bibliographystyle{siam}
\bibliography{/Users/contat/Dropbox/Articles/biblio.bib}
\end{document}